\documentclass[12pt]{amsart}
\usepackage{amsfonts,amsmath,amssymb}
\usepackage{mathabx}
\usepackage{mathrsfs}
\usepackage[latin1]{inputenc}
\usepackage[usenames,dvipsnames]{pstricks}
\newtheorem{theorem}{Theorem}[section]

\newtheorem{remark}[theorem]{Remark}

\newtheorem{Lemm}[theorem]{Lemma}

\renewcommand{\theequation}{\thesection.\arabic{equation}}

\newcommand{\p}{\partial}
\newcommand{\R}{\mathbb{R}}

\newcommand{\N}{\mathbb{N}}

\newcommand{\s}{\mathbb{S}}

\usepackage{graphicx}

\usepackage{mathrsfs}
\usepackage[T1]{fontenc}
\usepackage{color,xspace}
\usepackage[usenames,dvipsnames]{pstricks}
\usepackage{times}
\allowdisplaybreaks
\definecolor{rot}{rgb}{1,0,0}
\definecolor{hw}{rgb}{0,0,1}

\begin{document}
\renewcommand{\theequation}{\arabic{section}.\arabic{equation}}
\title[ ]
{ Increasing stability for inverse source problem with limited-aperture far field data at multi-frequencies }\author[Ben A\"icha, Hu and Si ]{ Ibtissem Ben A\"icha, Guang Hui Hu and Su Liang Si}

\address{  Ibtissem Ben A\"icha,  Higher Institute for Preparatory Studies in Biology-Geology (ISEP-BG), University of Carthage, Tunis, La Soukra 2036, Tunisia 
\& LAMSIN, National Engineering School of Tunis, B.P. 37, 1002 Tunis, Tunisia. 
	\newline
	\indent E-mail:{\tt \ ibtissem.benaicha@enit.utm.tn}}
\address{ Guang Hui Hu, School of Mathematical Sciences, Nankai University, Tianjin, 300071, China.
\newline
	\indent E-mail:{\tt \ ghhu@nankai.edu.cn}}
\address{ Su Liang Si, School of Mathematics and statistics, Shandong University of Technology, Shangdong, 255049, China.
\newline
	\indent E-mail:{\tt \ sisuliang@amss.ac.cn}}	

\maketitle
\begin{abstract}
We study the increasing stability of an inverse source problem for the Helmholtz equation from limited-aperture far field data at multiple wave numbers.  The measurement data are given by  the far field patterns $u^{\infty}(\hat{x},k)$ for all observation directions $\hat{x}$ in some neighborhood of a fixed direction $\hat{x}_0$ and for all wave numbers $k$ belonging to a  finite interval $(0,K)$.  In this paper, we  discuss  the increasing stability with respect to the width of the wavenumber interval $K>1$. In three dimensions we  establish  stability estimates of the $L^2$-norm and $H^{-1}$-norm of the source function from the far field data. The ill-posedness of the inverse source problem turns out to be of  H\"older type while increasing the wavenumber band $K$. We also discuss  an analytic continuation argument of  the far-field data with respect to the wavenumbers at a fixed direction. 


\end{abstract}

\section{Introduction}
\noindent In this paper, we deal with an inverse time-harmonic source  problem for the Helmholtz equation. Our aim is to determine the source function $f$ appearing in the following  equation defined in $\R^3$:
\begin{equation}\label{1}
\left\{
   \begin{array}{ll}
     (\Delta+k^2)u=-f & \mbox{in}\,\,\R^3, \\
\\  \displaystyle \lim r(\p_r u-iku)=0 & \mbox{as}\,\,r:=|x|\rightarrow +\infty.
   \end{array}
 \right.
\end{equation}
Here $k>0$ is the wave number, $u$ is the radiated wave field and $f\in L^{2}(\R^3)$ is a compactly supported function. We assume that the supp $(f)$ is contained in the ball $B(0,R)$ defined by 
$$B(0,R):=\{x\in\R^3: |x|<R\}$$
for some $R\geq1$. The condition at infinity in  (\ref{1}) is well known as the Sommerfeld radiation condition, leading to the asymptotic behavior 
\[
u(x)=\frac{e^{ik |x|}}{4\pi |x|}\left\{u^\infty(\hat{x}, k)+O\Big( \frac{1}{|x|} \Big) \color{black}\right\}\qquad\mbox{as}\quad |x|\rightarrow\infty,
\]
which holds uniformly in all directions $\hat{x}\in\mathbb{S}^2:=\{x\in \R^3: |x|=1\}$.
The function $\hat{x}\mapsto u^\infty(\hat{x}, k)$ is usually referred to as the far-field pattern or scattering amplitude of $u$ and the vector $\hat{x}\in\R^3$ is the observation direction. It is well known that $\hat{x}\mapsto u^\infty(\hat{x}, k)$ is a real-analytic function on $\mathbb{S}^2$. There exists  a unique solution $u\in H^2_{loc}(\R^3)$  to (\ref{1}), taking the explicit form
\begin{equation}\label{u}
u(x)=\int_{\R^3}G(x,y,k)\,f(y)\,dy,\quad x\in\R^3,
\end{equation}
where $G(x,y,k)$ denotes the free-space Green function of the Helmholtz equation and is given by 
$$G(x,y,k)=\frac{1}{4\pi} \frac{e^{ik|x-y|}}{|x-y|},\qquad x\neq y,\quad x, y\in \R^3.$$ 
From 
\[|x-y|=\sqrt{|x|^2-2x\cdot y+|y|^2}=|x|-\hat{x}\cdot y+O\Big(\frac{1}{|x|}\Big)\color{black}, \quad |x|\rightarrow\infty ,    \] 
we derive
\begin{equation}\label{e}
\frac{e^{ik|x-y|}}{|x-y|}=\frac{e^{ik|x|}}{|x|}\Big\{ e^{-ik \hat{x}\cdot y}+O\Big( \frac{1}{|x|}\Big)\color{black} \Big\} , \quad  |x|\rightarrow\infty.
\end{equation}
By \eqref{u} and \eqref{e}, 
the far-field pattern of $u$ is given by
\begin{equation}\label{2}
u^{\infty}(\hat{x},k)=\int_{\R^3} e^{-ik \hat{x}\cdot y} f(y)\,dy.
\end{equation}
Denote the Fourier transform of $f\in L^2(\R^3)$ by 
\begin{equation}
\hat{f}(\xi):=(2\pi)^{-\frac{3}{2}}\int_{\R^3}f(x)e^{-i\xi\cdot x}dx,\quad \xi\in \R^3.
\end{equation}
Thus, we get 
\[u^{\infty}(\hat{x},k)=(2\pi)^{\frac{3}{2}}\hat{f}(k\hat{x}).\]
For a small parameter $\delta>0$, let us denote by $\mathcal{V}_{\delta}(\hat{x}_0)$ a neighborhood  of a fixed observation direction $\hat{x}_0\in \mathbb{S}^2$, defined  as follows:
$$\mathcal{V}_{\delta}(\hat{x}_0):=\{\hat{x}\in\mathbb{S}^2: \,|\hat{x}-\hat{x}_0|<\delta\}.$$ 
Supposing for some $K>1$ that 
\[u^{\infty}(\hat{x},k)=0, \quad\mbox{for all}\quad \hat{x}\in\mathcal{V}_\delta(\hat{x}_0),\ k\in(0,K),\]
we deduce from the analyticity of $\hat{f}$ that
\[f(x)=0 \quad\mbox{for all}\quad x\in\R^3.\]
That is, the far field pattern $u^\infty(\hat{x},k)$ with $\hat{x}\in\mathcal{V}_\delta(\hat{x}_0)$ and  $k\in(0,K)$ uniquely determines the source term $f(x)$,  for  $x\in\R^3$.
The goal of this paper is to consider  the stability issue for the inverse problem of determining the source function $f$  from  limited-aperture far field data at multi-frequencies. 
More precisely, we are interested in treating the stability  issue with respect to the wavenumber/frequency interval $(0, K)$. Particularly, we are interested in the so-called increasing stability with partial data, showing that 
the ill-posedness of the inverse problem decreases as the width of the frequency band becomes larger.

We state the stability issue as follows:

\textit{\textbf{IP:}} Stably determine the source term $f$ appearing in (\ref{1}) from limited-aperture and multi-frequency far-field pattern data $u^\infty(\hat{x},k)$ for all $\hat{x}\in\mathcal{V}_{\delta}(\hat{x}_0)$ of a  fixed direction $\hat{x}_0\in\mathbb{S}^2$ and for all  wavenumbers $k\in (0,K)$ with some $K>1$. 

The overall idea used in resolving this problem is  mainly based on  an analytic continuation argument enabling the extension of the data  $u^{\infty}(\hat{x}, k)$ from $k\in (0,K)$ onto the real axis and from $\hat{x}\in \mathcal{V}_{\delta}(\hat{x}_0)$ onto the whole unit ball $\mathbb{S}^2$.   Indeed, if the available measurement data are confined to a finite frequency band and limited observation directions, one needs to extend them to larger intervals in order to be able to deal with theoretical and also numerical approaches.

The study of inverse coefficients and source  problems for partial differential equations is one of the most rapidly growing mathematical research area in the recent years,  and inverse source problems attracted recently much attention. Inverse source problems have indeed many applications in many fields, like for example in antenna synthesis \cite{[1]}, biomedical and medical imaging \cite{[2]} and  tomography \cite{[3],[4]}.

We state in brief some of the existing results that are relevant to the problem under investigation in this paper. It is well known  that there is an obstruction to uniqueness for inverse source problems for Helmholtz equations with a single frequency data. For the convenience of the reader, we refer to \cite{[6]} [Chapter 4] and \cite{BC}. However,  by considering multi-frequency measurements, the uniqueness can be proved. For this, one can see for example the recent works \cite{[5], [12]} in which uniqueness and stability results have been proved for the recovery of the source term from  knowledge of multi-frequency boundary measurements. In \cite{[7]}, the authors 
treated an interior inverse source problem for the Helmholtz equation from boundary Cauchy data for multiple wave numbers and they showed an increasing stability result for the problem under consideration. 
Interested readers can also see  \cite{[12]} that claims a uniqueness result and a numerical algorithm for recovering the location and the shape of a supported acoustic source function from boundary measurements at many frequencies.  See also \cite{[BLRX], [8],[11],[15],[16], [20]} and the references therein.
As for increasing stability results proved for coefficients inverse problems, we can refer for example to \cite{[Horst]} and \cite{[U]},  in which  inverse problems of recovering an electric potential appearing in a Schr\"odinger equation have been studied (see also the references therein). In \cite{[BT20]}, the increasing stability for the one-dimensional inverse medium  problem of recovering the refractive index is investigated. 
However, to the best of our knowledge, there is no analogue results for the Helmholtz equation using partial boundary measurement data or limited aperture far-field data. For this purpose, one needs  additionally to show the stability of the unique continuation of the measurement data  with respect to observation directions or wavenumbers.

Inspired by the works \cite{[BMBI], [7], [20]}, we address in the present paper, the increasing stability for  recovering the compactly supported function $f$ appearing in (\ref{1}) from the far-field pattern $u^{\infty}(\hat{x},k)$ for $\hat{x}\in \mathcal{V}_{\delta}(x_0)$ and $k\in(0,K)$. We show two stability estimates with respect to the frequency parameter $K$. These estimations consist of two parts: the first one is of H\"older type in terms of our data and the second one is  a logarithmic term that comes from  the high frequency tail of the function which decreases as $K$ increases, which makes the problem under investigation more stable.   
As a bi-product of our arguments,  we show a stability estimate in the unique continuation for the far field pattern with respect to $k$ or $\hat{x}$ but at fixed observation direction or fixed wavenumber respectively; see Theorem \ref{THM3} and Theorem \ref{THM4}.

The rest  of the paper is organized as follows.  In Section 2 we state our main results.   Sections 3 and 4 are devoted to the stability estimates of the $L^2$-norm and $H^{-1}$-norm of the source function. In Section 5, we present a proof to an analytic unique continuation argument of the far-field data with respect to the wavenumbers at a fixed observation direction or observation angles at a fixed wavenumber.

\section{Main results} 
Before stating our main results, we first introduce some notations.
Define a complex valued functional space
\[\mathcal{C}_{M,2n+1}=\left\{f\in H^{2n+1}(\R^3):\ \|f\|_{H^{2n+1}(\R^3)}\leq M,\ \mbox{supp} f\subset B_R \right\}\]
where $n\geq1$ is an integer and $M>1$ is a constant. 

For $K>1$,
let us  denote the $L^2$-norm of the limited-aperture far-field data by 
\begin{equation}\label{eps}
\epsilon:=\|u^{\infty}(\hat{x}, k)\|_{L^{\infty}(\mathcal{V}_\delta(\hat{x}_0)\times (0,K))}<1.
\end{equation}
Our first result concerns the stability estimate of the $L^2$-norm of $f$.
\begin{theorem}\label{THM1}
Let $f\in \mathcal{C}_{M,2n+1}$ and
let $u\in H^{2n+3}_{loc}(\R^{3})$ be the unique solution to the equation (\ref{1}). Then there exists $\alpha>0$ such that 
\begin{equation}\label{increasing2}
\|f\|^2_{L^2(\R^3)}\leq C\left[ K^3e^{2K(1-\alpha)}\epsilon^{2\alpha}+\frac{M^2}{(K^{\frac{2}{3}}|\ln\epsilon|^{\frac{1}{4}})^{4n-3}}\right],
\end{equation}
where $C>0$ is a constant independent of $K>0$.
\end{theorem}

Our second result provides an estimate of the $H^{-1}$-norm of $f$, where $f$ is only required to be a compactly supported $L^2$-function.
\begin{theorem}\label{THM2}
Let $f\in L^2(\R^3)$ be  compactly supported in $B_R$ with $M_0:=\|f\|_{L^2(\R^3)}$.
Let $u\in H^2_{loc}(\R^3)$ be the unique solution to the equation (\ref{1}). Then there exists $\alpha>0$ such that 
\begin{equation}\label{increasing1}
\|f\|^2_{H^{-1}(\R^3)}\leq C\left[  K^3e^{2K(1-\alpha)}\epsilon^{2\alpha}+\frac{M_0^2}{K^{\frac{4}{3}}|\ln\epsilon|^{\frac{1}{2}}}\right],
\end{equation}
where $C>0$ is a constant independent of $K>0$.
\end{theorem}
Comparing Theorems \ref{THM1} and \ref{THM1}, we require $f$ to be of high regularities in order to get the $L^2$-estimate of $f$. The inverse source problem appears to be more stable for smooth source functions. Theorems \ref{THM1} and \ref{THM1} show that
the ill-posedness of the inverse source problem decreases 
when  the wavenumber/frequency band $K$  increases. The inverse problem behaves like a H\"older  stability
using a large interval of wavenubmers, because the unstable logarithmic terms on the right hand sides of (\ref{increasing1}) and (\ref{increasing2}) decay to zero as $K$ grows. 
The exponential terms $e^{2K(1-\alpha)}$ in the H\"older part of the estimates (\ref{increasing1}) and (\ref{increasing2}) are due to the use of the analytic continuation argument and they do not pose any problem, since  $K$ is fixed in practice and so is the constant $e^{2K(1-\alpha)}$ (see also  \cite{[U]}). One should notice that these stability estimates imply the following uniqueness result with phased measurement data: the limited far-field data $u^{\infty}(\hat{x}, k)$ for all $\hat{x}\in \mathcal{V}_\delta(\hat{x}_0)$ and $k\in(0, K)$ uniquely determine $f$. 

As a by-product of the proofs of the aforementioned stability estimates, we get a quantitative unique continuation result for the function $k\mapsto u^\infty(\hat{x}_0, k)$ at a fixed $\hat{x}_0\in  \mathbb{S}^2$.

\begin{theorem}\label{THM3}
Let $\hat{x}_0\in\mathbb{S}^2$ be a fixed observation direction and suppose $a\geq \max\{2K, 1\}$. Then there exist  positive constants $C$ and $\gamma\in(0,1)$  such that 
\begin{equation}\label{proposition}
\|u^{\infty}(\hat{x}_0,\cdot)\|_{L^\infty(-a,a)}\leq C \,e^{a\,(1-\gamma)}\|u^{\infty}(\hat{x}_0,\cdot)\|_{L^\infty(0,K)}^\gamma,
\end{equation}
where $C$ and $\gamma$ depend only $K$ and $a$.
Moreover, the parameter $\gamma$  increases as $K/a$ increases.
\end{theorem}
Theorem \ref{THM3} suggests that,  if we only know the far field data at a fixed direction $\hat{x}_0\in\mathbb{S}^2$ and over some finite and small frequency interval $(0,K)$, then we can analytically extend the data with respect to $k$ in a larger domain. Moreover, the analytical extension becomes more stable if the ratio $K/a$ increases. 

Noting that for fixed $k>0$,
\begin{equation}
\lim\limits_{n \to \infty}\frac{(2k)^n}{n!}=0,
\end{equation}
we always find a $n_0\in\N$ such that 
\begin{equation}\label{n0}
\sup_{n\in\N}\frac{(2k)^n}{n!}=\frac{(2k)^{n_0}}{n_0!}.
\end{equation}

In the following Theorem, we deduce that the values of $\hat{f}(x)$ on $\{k\hat{x}: \ \hat{x}\in\mathcal{V}_{\delta_1}(\hat{x}_0)\}   \subset \partial B(0,k)$ with $\delta_1>\delta>0$ can be controlled by $\hat{f}(x)$ on $\{k\hat{x}: \ \hat{x}\in\mathcal{V}_{\delta}(\hat{x}_0)\}   \subset \partial B(0,k)$.
\begin{theorem}\label{THM4}
Let $k>0$ be a fixed parameter and let $\delta_1>\delta>0$. There exists $\gamma_1\in(0,1)$  such that 
\begin{equation}\label{M1}
\|u^{\infty}(\cdot,k)\|_{L^\infty(\mathcal{V}_{\delta_1}(\hat{x}_0))}\leq (2M_2)^{1-\gamma_1}\|u^{\infty}(\cdot,k)\|^{\gamma_1}_{L^\infty(\mathcal{V}_{\delta}(\hat{x}_0))}.
\end{equation}
Here \begin{equation}\label{M2}
M_2=\max\left\{\frac{(2k)^{n_0}}{n_0!}, 1\right\}\;\|f\|_{L^1(\R^3)}.
\end{equation}
Moreover, the parameter $\gamma$  increases as $\delta/\delta_1$ increases.
\end{theorem}
Theorem \ref{THM4} suggests that, if we only know the far field for a fixed $k>0$ and only in some neighborhood $\mathcal{V}_{\delta}(\hat{x}_0)$ of a fixed direction $\hat{x}_0$, then one can analytically extend it with respect to $\hat{x}$ to a larger domain of $\s^2$.

\section{Stability estimate of  the $L^2$-norm: proof of Theorem \ref{THM1}} 
This part is devoted to the stability estimate of  the $L^2$-norm of the source function
appearing in the equation (\ref{1}) from  knowledge of the far field pattern $u^{\infty}(\hat{x},k)$ given by (\ref{2}), for all $\hat{x}\in\mathcal{V}_\delta(\hat{x}_0)$  and all $k\in (0,K)$ for some $K>1$.  For this purpose we will use an analytic continuation argument given in \cite{[BMBI],[Vessela]}. We will also take inspirations from  \cite{[7]} and use their methods based  on  explicit bounds of the harmonic measure of $(0, K)$ in a sector of the complex plane $k=k_1+ik_2\in \mathbb{C}$.
\subsection{Preliminaries}\quad\\
The statements given in this subsection are crucial for the proof of our increasing stability estimate. Let us first introduce the following analytic continuation  statement that will be used thereafter. Set $B(0,c)=\{x\in \R^3: |x|<c\}$ for any $c>0$.
\begin{Lemm}\label{lem3.1}(See \cite{[BMBI]})
Let $\mathcal{O}$ be a non-empty  open set of the unit ball $B(0,1)\subset
\R^{3}$ and let $F$ be an analytic function in $B(0,2),$ that
satisfies
$$\|\p^{\gamma}F\|_{L^{\infty}(B(0,2))}\leq M_1\,|\gamma|!\,{\eta^{-|\gamma|}},\,\,\,\,\forall\,\gamma\in(\mathbb{N}\cup\{0\})^{3},$$
for some  $M_1>0$ and $\eta>0$ and $N(\eta)$. Then, we have 
$$\|F\|_{L^{\infty}(B(0,1))}\leq N\, M_1^{1-\alpha}\|F\|_{L^{\infty}(\mathcal{O})}^{\alpha},$$
where $\alpha\in(0,1)$ depends only  $\eta$ and $|\mathcal{O}|$.
\end{Lemm}
Using the above Lemma, we can control the Fourier transform of the source term by the far field pattern data in suitable norms as follows.
\begin{Lemm}\label{lem3.2}
There exists $\alpha(R, \delta)>0$ such that 
\begin{equation}\label{f}
\|\hat{f}\|_{L^{\infty}(B(0,K))}\leq C \,e^{K(1-\alpha)}\,\epsilon^{\alpha},
\end{equation}
where $C>0$ depends on $R$, $M_0$ and $\epsilon$ is defined by \eqref{eps}.
\end{Lemm}
\begin{proof}
From the definition of the far field pattern given by  (\ref{2}),  we have 
\begin{equation}\label{farfield pattern}
u^{\infty}(\hat{x},k)=\int_{\R^3} e^{-i k\hat{x}\cdot y}\,f(y)\,dy=(2\pi)^{\frac{3}{2}}\,\hat{f}(k\hat{x}),\quad \hat{x}\in \mathbb{S}^2,\,\,\,k\geq 0.
\end{equation}
Let us define the following set 
$$E_K:=\left\{\xi\in\R^3:\,\,\,\xi=k\hat{x},\,\, k\in[0,K),\,\, \hat{x}\in \mathcal{V}_{\delta}(\hat{x}_0)\right\}\subset \R^3.$$ \color{black}
In light of (\ref{farfield pattern}), we have 
\begin{equation}\label{4}
|\widehat{f}(\xi)|\leq C\, \|u^{\infty}\|_{L^{\infty}(\mathcal{V}_\delta(x_0)\times(0,K))}\quad \mbox{for\,\,all}\,\,\xi\in E_K,
\end{equation}
where $C>0$ depends on $R$.
In order to complete the proof of the Lemma, we need to extend the estimate (\ref{4}) that holds only in $E_K$ into  the ball $B(0,K)$. For this purpose let $c\in(0,1/2)$ and  let us consider the set $E_c\subset \R^3$.
Obviously, $KE_c=E_{cK} \subset E_K.$
For any $\xi\in\R^3$, we introduce the function 
\begin{equation}\label{Fr}
F_{K}(\xi)=\widehat{f}(K\xi),\qquad \xi\in \R^3.
\end{equation} Since the function $f$ is compactly supported,  the function $F_{K}$ is real analytic and it satisfies 
$$|\p_\xi^\gamma F_K(\xi)|\leq \Big|  \int_{B_R} (-i)^{|\gamma|} K^{|\gamma|} y^{|\gamma|} f(y)\,dy\Big|\leq ||f||_{L^1(B_R)}\,R^{|\gamma|} \,K^{|\gamma|}\leq ||f||_{L^1(B_R)}  e^{K}\,R^{|\gamma|}\,|\gamma|!.$$
where $\gamma=(\gamma_1, \gamma_2, \gamma_3)\in (\N\cup\{0\})^3$, $y^{|\gamma|}:=y_1^{\gamma_1}y_2^{\gamma_2} y_3^{\gamma_3}$, $|\gamma|=\gamma_1+\gamma_2+\gamma_3$ and we have used the inequality $K^m\leq e^K m !$ for all $m\in \N\cup\{0\}$. 
Applying Lemma \ref{lem3.1} to the function $F_K$ with  $M_1=||f||_{L^1(B_R)}\,e^{K}$,  $\eta=R^{-1}$ and with the set $\mathcal{O}$ given by
 $$\mathcal{O}:= B(0,1)\cap \displaystyle \displaystyle E_{c}\neq \emptyset,$$
there exist constants $C>0$ depending on $\eta$ and $\alpha=\alpha(\beta,\delta)\in (0,1)$  such that 
 $$\|F_K\|_{L^\infty (B(0,1))}\leq C \,e^{K(1-\alpha)} ||F_K||_{L^\infty (\mathcal{O})}^{\alpha}.$$
Thus, one gets in view of the identity (\ref{Fr}) the estimate
 $$\|\widehat{f}\|_{L^\infty(B(0,K))}\leq C e^{K(1-\alpha)}\|\widehat{f}\|_{L^\infty(K \mathcal{O})}.$$
Using the fact that  $$K \mathcal{O}=B(0,K)\cap K E_c \subset {E}_K,$$ \color{black}
we can deduce in light of (\ref{4}),
\begin{equation}\label{5}
\|\widehat{f}\|_{L^{\infty}(B(0,K))}\leq C\, e^{K(1-\alpha)} \|u^{\infty}\|^{\alpha }_{L^{\infty}(\mathcal{V}_{\delta}(\hat{x}_0)\times(0,K))}.
\end{equation}
This completes the proof of the Lemma.
\end{proof}

Now we prove the main statement  of the present manuscript by taking inspiration from \cite{[7]}.
Denote
\begin{equation}\label{I}
I_1(k)=\int_{|\xi|\leq k}|\hat{f}(\xi)|^2d\xi=\int_{0}^k \int_{\mathbb{S}^2}|\hat{f}(l\theta)|^2\,l^2\,d\theta \,dl. 
\end{equation} 
Since the integrand is an entire analytic function of $\xi$, the integral in (\ref{I}) with respect to $l$ can be taken over any path jointing the points $0$ and $k$ of the complex plane. Thus $I_1(k)$ is an entire function of $k=k_1+ik_2$ and the following elementary estimates hold.

\begin{Lemm}\label{lem3.3}
Let $k=k_1+i k_2\in\mathbb{C}$ and $f\in L^2(B_R)$. Then 
\begin{equation}\label{2.8}
\Big|I_1(k)\Big| \leq C \,|k|^2 e^{2R|k_2|}.
\end{equation}
Here the constant $C>0$ depends on $\|f\|_{L^2(\R^3)}$ and $R$.
\end{Lemm}
\begin{proof}
Using the change of variables $l=ks$ for $s\in (0,1)$ and H\"older inequality, one can easily derive 
\begin{eqnarray}\label{Ik}
\Big|\int_0^k \int_{\mathbb{S}^2} |\hat{f}(l\theta)|^2\,l^2\,d\theta \,dl\Big|&=&\Big|  \int_{0}^1\int_{\mathbb{S}^2}|\hat{f}(ks\theta)|^2 (ks)^2\,d\theta\,ds \Big|\cr
&\leq& C \int_{0}^1\int_{\mathbb{S}^2} |k|^2  \int_{B_R} |f(x)|^2 \, |e^{2\, k_2 s\, x\cdot\theta}|\,dx\,d\theta\,ds \\
&\leq& C\, ||f||_{L^2(\R^3)}^2\, |k|^2\, e^{2\,R\, |k_2| }.
\end{eqnarray}
Here we used the fact that  $|e^{2\, k_2 s\,x\cdot\theta}|\leq e^{2\,R\,|k_2|}$. Then one gets the desired estimate (\ref{2.8}). 
\end{proof}
\noindent Now we are in a position to state an important Lemma  given in \cite [Lemma 3.2]{[7]} and also in \cite[Page 59]{[6]}. This Lemma plays a crucial role in our proof. Let us first introduce this complex sector  $\mathcal{S}$ by
$$\mathcal{S}:=\Big\{k=k_1+ik_2\in\mathbb{C}:\  |\arg\, k|< \frac{\pi}{4}\Big\}.$$
\begin{Lemm}\label{lem3.4}(\cite{[6],[7]})
Let $J(k)$ be an analytic function in $\mathcal{S}$ and continuous in $\overline{\mathcal{S}}$ satisfying 
\begin{equation*}
\left\{
   \begin{array}{ll}
     |J(k)|\leq a, & k\in (0,K],\\
     \\

 |J(k)| \leq B , & k\in \mathcal{S},\\
 \\
 |J(0)|=0,
   \end{array}
 \right.
\end{equation*}
where $0<a<B$. 
Then there exists a function $\mu(k)$ (a harmonic measure) satisfying 
\begin{equation}\label{mu}
\left\{
   \begin{array}{ll}
     \mu(k)\geq \frac{1}{2}, & k\in (K,2^{1/4} K),\\
     \\
\mu(k)\geq \frac{1}{\pi}((\frac{k}{K})^4-1)^{-1/2} , & k\in (2^{1/4} K,+\infty),
   \end{array}
 \right.
\end{equation}
such that 
$$|J(k)|\leq B \,a^{\mu(k)},\quad \forall\,k\in (K,+\infty).$$
\end{Lemm}
\noindent For $k=k_1+i k_2\in \mathcal{S}$, let us denote by $J(k)$  the following integral function
$$J(k)=e^{-2(R+1)k}I_1(k)=e^{-2(R+1)k}\int_0^k\int_{\mathbb{S}^2}|\hat{f}(l\theta)|^2\,l^2\,d\theta \,dl.$$

Using Lemmas \ref{lem3.2} and \ref{lem3.4}, we show the following statement. 
\begin{Lemm} \label{lem3.5}
There exist $\alpha\in(0,1)$ and a function $\mu(x)$ satisfying (\ref{mu}) such that 
\begin{equation}\label{eq2.13}
|J(k)|\leq C\, e^{2(1-\alpha) K \mu(k)} \epsilon^{2\alpha \mu(k)},
\end{equation}
holds true for any $k>K$. Here $C\!\!=\!\!C(M_0, R)>0$ with $M_0=\|f\|_{L^2(\R^3)}$.
\end{Lemm}
\begin{proof}
Since $|k_2|\leq k_1$ for any $k\in \mathcal{S}$, then in view of Lemma \ref{lem3.3}, we have 
\begin{equation}\label{est1}
|J(k)|\leq C_2|k|^2 e^{2 R |k_2|}|e^{-2 (R+1)k}|\leq  C_2, 
\end{equation}
where $C_2>0$ depending on $M_0$ and $R$.
On the other hand, in light of Lemma \ref{lem3.2} and using the fact that  $k^3\leq e^{2(R+1)k}$ and $k>0$, one has for any $k\in (0,K]$, 
\begin{equation}\label{est2}
|J(k)|\leq \frac{4\pi}{3}e^{-2 (R+1) k}k^3\|\hat{f}\|^2_{L^{\infty}(B(0,K))}\leq C_1\, e^{2 K(1-\alpha)} \epsilon^{2\alpha}
\end{equation}
where $C_1>0$ depends on $R$. 
Let us denote by $C=\max\{C_1, C_2\}$ depending on $M_0$ and $R$. 
Therefore, from the estimates (\ref{est1}) and  (\ref{est2}) and Lemma \ref{lem3.4} with $a=C e^{2 K(1-\alpha)} \epsilon^{2\alpha}$ and  $B=C$, we know that there exists a function $\mu(k)$ satisfying (\ref{mu}) such that (\ref{eq2.13}) holds true.
\end{proof}
Denote
\[I_2(s)=\int_{|\xi|>s}|\hat{f}(\xi)|^2d\xi=\int_s^\infty\int_{\s^2}|\hat{f}(k\hat{x})|^2k^2\;d\hat{x}dk.\] 
Now we estimate $I_2(s)$.
\begin{Lemm}\label{lem3.8}
Let $f\in \mathcal{C}_{M,2n}$.
For any $s\geq1$, we have 
\begin{equation}\label{I_2}
|I_2(s)|\leq\frac{C}{s^{4n-3}},
\end{equation}
where $C>0$ depends on $M$ and $n$.
\end{Lemm}
\begin{proof} 
Let $\Delta f(y)=\partial_{y_1}^2f(y)+\partial_{y_2}^2f(y)+\partial_{y_3}^2f(y)$.
Denote $\Delta^nf(y):=\underbrace{\Delta\cdot\cdot\cdot\Delta(\Delta }_{n}f(y))$.
Since
$\hat{f}$ is the Fourier transform of $f$ given by 
\begin{equation}\label{Ff}
\hat{f}(k\hat{x})=(2\pi)^{-\frac{3}{2}}\int_{\R^3}e^{-ik\hat{x}\cdot y}\; f(y)dy.
\end{equation}
Multiplying $(-ik\hat{x}_j)^{2}$,\ $j=1, 2, 3$ on both sides of (\ref{Ff}), and adding the three equations now gives 
\begin{eqnarray}\label{3i}
\big((-ik\hat{x}_1)^{2}+(-ik\hat{x}_2)^{2}&+&(-ik\hat{x}_3)^{2}\big)\hat{f}(k\hat{x})\cr
&=&\big((-ik\hat{x}_1)^{2}+(-ik\hat{x}_2)^{2}+(-ik\hat{x}_3)^{2}\big)(2\pi)^{-\frac{3}{2}}\int_{\R^3}e^{-ik\hat{x}\cdot y}\; f(y)dy\cr
&=&(2\pi)^{-\frac{3}{2}}\int_{\R^3}\Delta e^{-ik\hat{x}\cdot y}\; f(y)dy\cr
&=&(2\pi)^{-\frac{3}{2}}\int_{\R^3}(e^{-ik\hat{x}\cdot y})\; \Delta f(y)dy.
\end{eqnarray}
Since $\hat{x}_1^2+\hat{x}_2^2+\hat{x}_3^2=1$, we know from (\ref{3i}) that
\begin{equation}
(-ik)^{2}\hat{f}(k\hat{x})=(2\pi)^{-\frac{3}{2}}\int_{\R^3}(e^{-ik\hat{x}\cdot y})\;\Delta f(y)dy.
\end{equation}
Similarly, repeating the above process $n$ times shows that
\begin{equation}\label{nf}
(-ik)^{2n}\hat{f}(k\hat{x})=(2\pi)^{-\frac{3}{2}}\int_{\R^3}(e^{-ik\hat{x}\cdot y})\;\Delta^n f(y)dy.
\end{equation} 
Hence one immediately has from (\ref{nf}) that 
\begin{equation}
|\hat{f}(k\hat{x})|\leq\frac{C_3}{k^{2n}}\quad\mbox{for\ all}\quad \hat{x}\in \s^2, \ k>0,
\end{equation}
where $C_3>0$ depends on $M$ and $R$. This leads to
\begin{eqnarray}
I_2(s)&=&\int_s^\infty\int_{\s^2}|\hat{f}(k\hat{x})|^2k^2ds(\hat{x})dk\cr
&\leq &4\pi{C_3}^2\int_s^\infty\frac{1}{k^{4n-2}}dk\cr
&\leq &\frac{4\pi{C_3}^2}{4n-3}\frac{1}{s^{4n-3}}.
\end{eqnarray}
Setting $C=\frac{4\pi{C_3}^2}{4n-3}$, we obtain (\ref{I_2}).
\end{proof}
\subsection{Proof of Theorem \ref{THM1}}\quad \\ 
Set $\alpha=\alpha(R, \delta)>0$, which is defined in Lemma \ref{lem3.2}.

We assume that $\epsilon<e^{-1}$, otherwise the estimate is obvious. Let 
\begin{eqnarray}
s=
\begin{cases}
\frac{\alpha^{\frac{1}{3}}}{(4\pi(R+1))^{\frac{1}{3}}}K^{\frac{2}{3}}|\ln\epsilon|^{\frac{1}{4}} & \textrm{if}\ |\ln\epsilon|^{\frac{1}{4}}>2^{\frac{1}{4}}\frac{K^{\frac{1}{3}}(4\pi(R+1))^{\frac{1}{3}}}{\alpha^{\frac{1}{3}}},\\
K & \textrm{if} \ |\ln\epsilon|^{\frac{1}{4}}\leq2^{\frac{1}{4}}\frac{K^{\frac{1}{3}}(4\pi(R+1))^{\frac{1}{3}}}{\alpha^{\frac{1}{3}}}.
\end{cases}
\end{eqnarray}

Case (i): $|\ln\epsilon|^{\frac{1}{4}}>2\frac{K^{\frac{1}{3}}(4\pi(R+1))^{\frac{1}{3}}}{\alpha^{\frac{1}{3}}}$. 
One can check that 
\[2\alpha|\ln\epsilon|-2K(1-\alpha)>\alpha|\ln\epsilon|.\]
Thus, using Lemma \ref{lem3.5} we obtain
\begin{eqnarray}
|I_1(s)|&\leq & C\, e^{2(R+1)s} e^{2(1-\alpha) K \mu(k)}\epsilon^{2\alpha \mu(k)}\cr
&\leq & Ce^{-\big((2\alpha|\ln\epsilon|-2K(1-\alpha))\mu(s)-2(R+1)s\big)}\cr
&\leq & Ce^{-\big((2\alpha|\ln\epsilon|-2K(1-\alpha))\frac{1}{\pi}(\frac{K}{s})^2-2(R+1)s\big)}\cr
&\leq & Ce^{-\big(\frac{\alpha^{\frac{1}{3}}}{\pi}(4\pi(R+1))^{\frac{2}{3}}K^{\frac{2}{3}}|\ln\epsilon|^{\frac{1}{2}}(1-|\ln\epsilon|^{-\frac{1}{4}})\big)},
\end{eqnarray}
where $C$ is defined in Lemma \ref{lem3.5}. 

Noting that 
$\frac{1}{2}|\ln\epsilon|^{-\frac{1}{4}}<\frac{1}{2}$,
we have 
\[|I_1(s)|\leq Ce^{-\big(\frac{\alpha^{\frac{1}{3}}(4\pi(R+1))^{\frac{2}{3}}}{2\pi}K^{\frac{2}{3}}|\ln\epsilon|^{\frac{1}{2}}\big)}.\]
Using the elementary inequality 
\[e^{-t}\leq \frac{(12n-9)!}{t^{3(4n-3)}}, \quad t>0,\]
we get 
\begin{equation}
|I_1(s)|\leq\frac{C}{\big(K^2|\ln\epsilon|^{\frac{3}{2}}\big)^{4n-3}},
\end{equation}
where $C>0$ depends on $M_0$, $R$, $n$ and $\alpha$. Since $K^{\frac{2}{3}}|\ln\epsilon|^{\frac{1}{4}}\leq K^2|\ln\epsilon|^{\frac{2}{3}}$ when $K>1$ and $|\ln\epsilon|>1$, we have 
\begin{equation}
|I_1(s)|\leq\frac{C}{\big(K^{\frac{2}{3}}|\ln\epsilon|^{\frac{1}{4}}\big)^{4n-3}}.
\end{equation}
Case (ii): $|\ln\epsilon|^\frac{1}{4}\leq2\frac{K^{\frac{1}{3}}(4\pi(R+1))^{\frac{1}{3}}}{\alpha^{\frac{1}{3}}}$. In this case $s=K$. We have from (\ref{f}) that 
\[|I_1(s)|=|I_1(K)| \leq CK^3e^{2K(1-\alpha)}\epsilon^{2\alpha}.\]

Combining the estimates of $I_1(s)$ and $I_2(s)$ and applying Lemma \ref{lem3.8}, we obtain 
\begin{eqnarray}\label{Hest0}
\|f\|^2_{L^2(\R^3)}&=&\|\widehat{f}\|^2_{L^2(\R^3)}\cr
&=&\int_{|\xi|\leq s}|\widehat{f}(\xi)|^2\,d\xi+\int_{|\xi|>s}|\widehat{f}(\xi)|^2\,d\xi\cr
&=&I_1(s)+I_2(s)\cr
&\leq& CK^3e^{2K(1-\alpha)}\epsilon^{2\alpha}+\frac{C}{(K^{\frac{2}{3}}|\ln\epsilon|^{\frac{1}{4}})^{4n-3}}+\frac{M^2}{\big(2^{-\frac{1}{4}}\alpha^{-\frac{1}{3}}(4\pi(R+1))^{-\frac{1}{3}}K^{\frac{2}{3}}|\ln\epsilon|^{\frac{1}{4}}\big)^{4n-3}}\cr
&\leq& C\left[ K^3e^{2K(1-\alpha)}\epsilon^{2\alpha}+\frac{M^2}{(K^{\frac{2}{3}}|\ln\epsilon|^{\frac{1}{4}})^{4n-3}}\right].
\end{eqnarray}
Thus we obtain the stability estimate (\ref{increasing2}).
This completes the proof.

\section{Stability estimate of the $H^{-1}$-norm: proof of Theorem \ref{THM2}}\quad \\
In this section we suppose that $f\in L^2(\R^2)$ is compactly supported in $B(0, R)$. Set $M_0=\|f\|_{L^2(\R^3)}$.
Let us first decompose the $H^{-1}(\R^{3})$ norm of $f$ into the following way
\begin{eqnarray}\label{31}
\|f\|^2_{H^{-1}(\R^{3})}&=&\int_{|\xi|\leq s}(1+|\xi|^2)^{-1} |\widehat{f}(\xi)|^2\,d\xi+\int_{|\xi|>s} (1+|\xi|^2)^{-1}|\widehat{f}(\xi)|^2\,d\xi.
\end{eqnarray}
We start by examining the last integral. The Parseval-Plancherel theorem and the H\"older inequality imply 
\begin{eqnarray}\label{32}
\int_{|\xi|>s}(1+|\xi|^2)^{-1}|\widehat{f}(\xi)|^2\,d\xi &\leq & \frac{1}{s^2}\int_{|\xi|>s}|\widehat{f}(\xi)|^2\,d\xi
\leq \frac{1}{s^2}\int_{\R^3}|\widehat{f}(\xi)|^2\,d\xi \\ \nonumber
&=&\frac{1}{s^2}\int_{\R^3}|f(x)|^2\,dx=\frac{M^2_0}{s^2}.
\end{eqnarray}
Further, in light of Lemma \ref{lem3.2}, if $0<s\leq K$, we have
\begin{equation}\label{33}
\int_{|\xi|\leq s}(1+|\xi|^2)^{-1} |\widehat{f}(\xi)|^2\,d\xi\leq\int_{|\xi|\leq K}|\widehat{f}(\xi)|^2\,d\xi\leq CK^3e^{2K(1-\alpha)}\epsilon^{2\alpha}.
\end{equation}
If $s>K$, in a similar way to Lemma \ref{lem3.5}, we get 
\begin{eqnarray}\label{34}
I_1(s)=\int_{|\xi|\leq s}(1+|\xi|^2)^{-1} |\widehat{f}(\xi)|^2\,d\xi&\leq & C\, e^{2(R+1)s} e^{2(1-\alpha) K \mu(k)}\epsilon^{2\alpha \mu(k)}\cr
&\leq & Ce^{-\big((2\alpha|\ln\epsilon|-2K(1-\alpha))\mu(s)-2(R+1)s\big)}\cr
&\leq & Ce^{-\big(\frac{\alpha^{\frac{1}{3}}(4\pi(R+1))^{\frac{2}{3}}}{2\pi}K^{\frac{2}{3}}|\ln\epsilon|^{\frac{1}{2}}\big)}.
\end{eqnarray}
Using the elementary inequality 
\[e^{-t}\leq \frac{6!}{t^{6}}, \quad t>0,\]
we get 
\begin{equation}\label{35}
|I_1(s)|\leq\frac{C}{\big(K^2|\ln\epsilon|^{\frac{3}{2}}\big)^2},
\end{equation}
where $C>0$ depends on $R$, $\alpha$ and $M_0$.
Since $K^{\frac{2}{3}}|\ln\epsilon|^{\frac{1}{4}}\leq K^2|\ln\epsilon|^{\frac{2}{3}}$ when $K>1$ and $|\ln\epsilon|>1$, we obtain 
\begin{equation}\label{36}
|I_1(s)|\leq\frac{C}{\big(K^{\frac{2}{3}}|\ln\epsilon|^{\frac{1}{4}}\big)^2}.
\end{equation}
Combining the proof of Theorem \ref{THM1} with (\ref{31})-(\ref{36}), one gets 
\begin{eqnarray}
\|f\|^2_{H^{-1}(\R^{3})}&=&\int_{|\xi|\leq s}(1+|\xi|^2)^{-1} |\widehat{f}(\xi)|^2\,d\xi+\int_{|\xi|>s} (1+|\xi|^2)^{-1}|\widehat{f}(\xi)|^2\,d\xi\cr
&\leq& CK^3e^{2K(1-\alpha)}\epsilon^{2\alpha}+\frac{C}{(K^{\frac{2}{3}}|\ln\epsilon|^{\frac{1}{4}})^{2}}+\frac{M_0^2}{\big(2^{-\frac{1}{4}}\alpha^{-\frac{1}{3}}(4\pi(R+1))^{-\frac{1}{3}}K^{\frac{2}{3}}|\ln\epsilon|^{\frac{1}{4}}\big)^{2}}\cr
&\leq& C\left[  K^3e^{2K(1-\alpha)}\epsilon^{2\alpha}+\frac{M_0^2}{K^{\frac{4}{3}}|\ln\epsilon|^{\frac{1}{2}}}\right].
\end{eqnarray}
Thus, we obtain the stability estimate (\ref{increasing1}).
This completes the proof. 

\section{Analytic continuation with respect to wavenumbers and observation angles: proof of Theorems \ref{THM3} and \ref{THM4}}
In this section,  we give an analytic continuation property for the far field data $u^{\infty}(\hat{x}_0,k)$ when  the direction $\hat{x}_0\in\mathbb{S}^2$ or $k>0$ is fixed. We recall that the far field pattern is defined in terms of the compactly supported function $f\in L^2(B_R)$,  given in (\ref{2}). We first show a stability estimate in the unique continuation for the far field pattern  with respect to $k$.  Let us first recall the following lemma proved in \cite{[BMBI], [Vessela]}.

\begin{Lemm}\label{LemmaA.2}
Let $\varphi$ be an analytic function in $[-1,1]$, and $I$ an open interval in
$[-1,1]$. We assume that there exist positive constants $M_3$ and $\rho$ such
that
\begin{equation}\label{A.49}
|\varphi^{(\kappa)}(s)|\leq \frac{M_3 \kappa!}{ (2\rho)^{\kappa}},\,\,\,\,\,\,\,\,\, \kappa\geq 0,\,\,s\in [-1,1].
\end{equation}
Then, there exist positive constants $N=N(\rho,|I|)$ and $\gamma=\gamma(\rho,|I|)$ such that 
\begin{equation}\label {A.50}
 |\varphi(s)|\leq N
\|\varphi\|_{L^{\infty}(I)}^{\gamma}M_3^{1-\gamma},\,\,\,\,\,\mbox{for\,\,all}\,\,\,\,s\in[-1,1].
\end{equation}
Here $N$ and $\gamma$ are given by 
\begin{eqnarray}\label{C}
N=2\delta_1^{-\gamma}(\ln L)^{\gamma}+(\ln\frac{8}{7})^{\gamma}(\ln L)^{-\gamma},\quad
\gamma:= (\ln\frac{8}{7}L)^{-1}\ln\frac{8}{7}.
\end{eqnarray} Here  $L:=6/|J_{s_{j_0}}|$ where the interval $J_{s_{j_0}}$ is defined as follows.
Let $$s_j=-1+\frac{(2j-1)\rho}{5}, \quad \frac{5}{\rho}\leq n_1\leq\frac{5}{\rho}+\frac{1}{2}, \quad I_j=[s_j-\frac{\rho}{5},s_j+\frac{\rho}{5}).$$
Choose $j_0\in(1,\cdot\cdot\cdot,n_1)$ such that $|I_{j_0}\cap I|=\max_{1\leq j\leq n_1}|I_{j_0}\cap I|$ and define $J_{s_{j_0}}=\frac{1}{\rho}(I_{j_0}\cap I-s_{j_0})$.
We consider the function $g$ defined as follows 
$$g(s)=\frac{\varphi(s_{j_0}+\rho s)}{2M_3}.$$
The parameter $\delta_1$ is required to satisfy 
\[
\ln\frac{8}{7}-e\|g\|_{L^\infty(J_{s_{j_0}})}\ln L\geq\delta_1>0.\]
\end{Lemm}
The above Lemma can be regarded as the one-dimensional analogue of Lemma \ref{lem3.6} with more explicit dependence of $N$ and $\gamma$ on $|I|$ and $\rho$.
Now, we present the proof of  Theorem \ref{THM3} by applying Lemma \ref{LemmaA.2}.


{\bf Proof of Theorem \ref{THM3}.} Let $\hat{x}_0\in \mathbb{S}^2$ be fixed and set $N_0:=a/K\geq 2$.
 For  a large parameter $1\leq a=N_0 K$, we define the following function
$$\varphi_a(k):=u^{\infty}(\hat{x}_0, k a),\quad k\in(-1,1).$$
One can see that $\varphi_a$ is analytic and  the following identity holds true
$$|\varphi_a^{(\kappa)}(k)|=|\p_k^{\kappa}\int_{\R^3} e^{-ika\hat{x}_0\cdot y} \,f(y)\,dy|\leq e^a \kappa!\, \|f\|_{L^1(B_R)}\, R^{\kappa}.$$
Applying the analytic continuation argument of Lemma \ref{LemmaA.2} with 
$$I:=\Big(0,\frac{1}{2N_0}\Big),\quad M_3=e^a\; \|f\|_{L^1(B_R)}\quad \mbox{and}\quad \rho=R^{-1},$$ there exist positive constants  $C$ and $\gamma$ given by (\ref{C}) and depending on $R$ and $N_0$ such that 
$$\|\varphi_a\|_{L^\infty(-1,1)}\leq  C e^{a(1-\gamma)}\,\|\varphi_a\|^{\gamma}_{L^\infty(I)}.$$
Now, using the fact that $a I\subset (0,K)$, one can see that 
$$\|u^{\infty}(\hat{x}_0,\cdot)\|_{L^\infty((-a,a)}\leq C e^{a\,(1-\gamma)}\,||u^\infty(\hat{x},\cdot)||_{L^\infty(0,K)}^\gamma.$$
Moreover, from Lemma \ref{LemmaA.2} we know 
\begin{equation}
\gamma:= (\ln\frac{48}{7|J_{s_{j_0}}|})^{-1}\ln\frac{8}{7},\quad 
J_{s_{j_0}}:=\frac{1}{\rho}(I_{j_0}\cap I-s_{j_0}). 
\end{equation}
where $s_{j_0}$ is defined in Lemma \ref{LemmaA.2}.
It is easy to know there exists a large constant $C_0$ such that $$C_0\frac{|I|}{\rho}\leq|J_{s_{j_0}}|\leq\frac{|I|}{\rho}.$$
Thus,  the exponent $\gamma$ satisfies
\begin{equation}\label{AI}
(\ln\frac{48\rho}{7C_0|I|})^{-1}\ln\frac{8}{7}\leq\gamma\leq(\ln\frac{48\rho}{7|I|})^{-1}\ln\frac{8}{7}.
\end{equation}
Since $|I|=\frac{1}{2N_0}$, the inequalities in (\ref{AI}) can be rewritten as 
\[(\ln\frac{96\rho N_0}{7C_0|I|})^{-1}\ln\frac{8}{7}\leq\gamma\leq(\ln\frac{96\rho N_0}{7})^{-1}\ln\frac{8}{7}.\]
This means that $\gamma$ is getting smaller (resp. bigger),
if $N_0$ becomes larger (resp. smaller). Equivalently, we conclude that 
 if $K$ increases or $a$ decreases , the exponent $\gamma$ will increase. 

The following Lemma is proved in \cite{[Vessela]}.
\begin{Lemm}\label{lem3.6}
Let $\Omega$ be a connected bounded, open set in $\R^n$ ($n=2,3$) such that for a positive number $r_0$ the set 
$$
\Omega_r:=\left\{x\in\Omega:\; d(x,\partial\Omega)>r\right\},\quad d(x,\partial\Omega):=\inf\limits_{y\in\partial\Omega}|x-y|
$$
is connected for every $r\in[0,r_0]$. Let $E\subset\Omega$ be an open set such that $d(E,\partial\Omega)\geq d_0>0$, and $f$ an analytic function on $\Omega$, having the following property 
\begin{equation}
|\partial^\alpha f(x)|\leq\frac{M_2\alpha!}{\rho^{|\alpha|}} \qquad \textrm{for all}\quad x\in\Omega, \ \alpha\in(\N\cup\{0\})^n,
\end{equation}
where $\rho$, $M_2$ are positive numbers.
Then
\begin{equation}
|f(x)|\leq (2M_2)^{1-\gamma_1}\big(\sup_{E}|f|\big)^{\gamma_1}\qquad\mbox{for all}\quad x\in \Omega,
\end{equation}
where diam $\Omega=\sup\{|x-y|:\ x,\ y\in\Omega\}$, $\gamma_1\in(0,1)$ depends only on  $\frac{|E|}{|\Omega|}$, $d_0$, diam $\Omega$, $n$, $r_0$, $\rho$ and $d(x,\partial\Omega)$, 
where $|E|$ and $|\Omega|$ denote the Lebesgue measure of $E$ and $\Omega$ respectively.
\end{Lemm}

\begin{remark}\label{rem}
By the proof of Lemma \ref{lem3.6} in \cite{[Vessela]}, we know that if $\frac{|E|}{|\Omega|}$ increases, the exponent $\gamma_1$ will increase.
\end{remark}
Now we show a stability estimate in the unique continuation for the far field pattern  with respect to observation angles $\hat{x}$.

{\bf Proof of Theorem \ref{THM4}.}
We use polar coordinates to represent $\mathcal{V}_{\delta}(\hat{x}_0)$ such that 
\[\mathcal{V}_{\delta}(\hat{x}_0)=\left\{(\cos\theta\cos\varphi,\sin\theta\cos\varphi,\sin\varphi):\ (\theta,\varphi)\in E\right\},\]
where
\[E=\{(\theta,\varphi):\ \theta_1<\theta<\theta_2,\ \varphi_1<\varphi<\varphi_2\}\]
for some $0<\theta_1<\theta_2<\pi$, $0<\varphi_1<\varphi_2<2\pi$.
Suppose $\Omega\supset E$ is an open set such that $d(E,\Omega)\geq d_0>0.$ Then
\[\mathcal{V}_{\delta}(\hat{x}_0)\subset\mathcal{V}_{\delta_1}(\hat{x}_0)=\{(\cos\theta\cos\varphi,\sin\theta\cos\varphi,\sin\varphi):\ (\theta,\varphi)\in \Omega\}.\]
Denote
\[
H(\hat{x})=H(\theta,\varphi)=\int_{\R^3}e^{-ik\hat{x}\cdot y}f(y)dy,\]
where 
$\hat{x}=(\cos\theta\cos\varphi,\sin\theta\cos\varphi,\sin\varphi).$
If $k>1$, since supp $f(y)\subset B_R$, we have for $\beta\in (\N\cup\{0\})^2$ that
\begin{eqnarray}\nonumber
\Big|\partial^\beta H(\theta,\varphi)\Big|&=&\Big|\int_{\R^3}\partial^\beta e^{-ik\hat{x}\cdot y}f(y)dy\Big|
\leq \int_{\R^3}|f(y)|dy\; (2k)^{|\beta|}(2R)^{|\beta|}\\ \nonumber
&\leq & \|f\|_{L^1(\R^3)} \frac{(2k)^{|\beta|}}{\beta!}(2R)^{|\beta|}\beta! \\ \label{p1}
&\leq & \|f\|_{L^1(\R^3)} \frac{(2k)^{|n_0|}}{n_0!}(2R)^{|\beta|}\beta! ,\quad
\end{eqnarray}
for all $(\theta,\varphi)\in \Omega$, where $n_0\in \N$ is given by \eqref{n0}.
If $0<k\leq1$,then 
\begin{eqnarray}\label{p2}
\Big|\partial^\beta H(\theta,\varphi)\Big|=\Big|\int_{\R^3}\partial^\beta e^{-ik\hat{x}\cdot y}f(y)dy\Big| 
\leq \|f\|_{L^1(\R^3)}\;(2R)^{|\beta|}\beta!\quad \mbox{for all}\ (\theta,\varphi)\in \Omega.
\end{eqnarray}
Set  $\rho=(2R)^{-1}$ and let $M_2$ be defined by \eqref{M2}.
Combining estimates (\ref{p1}) and (\ref{p2}) with Lemma \ref{lem3.6}, we have
\begin{equation}\label{p3}
|H(\theta,\varphi)|\leq (2M_2)^{1-\gamma_1}\big(\sup_{E}|H(\theta,\varphi)|\big)^{\gamma_1},
\end{equation}
for $(\theta,\varphi)\in \Omega$. 
In fact, according to Remark \ref{rem}, the parameter $\gamma$ increases as $\delta$ increases or  $\delta_1$ decreases.
Replacing $H(\theta,\varphi)$ by $u^{\infty}(\hat{x},k)$ in (\ref{p3}), we deduce (\ref{M1}).




\section*{Acknowledgment}
The work of G. Hu is partially supported by the National Natural Science Foundation of China (No. 12071236) and the Fundamental Research Funds for Central Universities in China (No. 63233071). The work of S. Si is supported by the Natural Science Foundation of Shandong Province, China (No. ZR202111240173). The first author would like to thank Professor Mourad Bellassoued for many helpful discussions.

\end{document}